\documentclass[12pt]{article}
\usepackage {amsfonts}
\usepackage[english]{babel}
\newtheorem{Th}{Theorem}
\newtheorem{Le}{Lemma}

\def\R{{\mathbb R}}
\def\C{{\mathbb C}}
\def\N{{\mathbb N}}
\def\Z{{\mathbb Z}}

\def\e{\varepsilon}

\def\<{\langle }
\def\>{\rangle }

\def\eq{eqnarray}
\def\l{\lambda}
\def\a{\alpha}

\def\card{{\rm card}}
\def\arg{{\rm arg}}
\def\E{{\rm E}}

\begin{document}
\title{
    \bf Zero sets of entire functions of exponential type
    with some conditions on real axis}

\author{Favorov S.Ju.}

\date{}

\maketitle

\begin{abstract}

We give a complete description of zero sets for some well-known
subclasses of entire functions of exponential growth.

\end{abstract}



\bigskip

By $B$ we denote the class of entire functions of exponential type
bounded on the real axis;

by $C$ denote the class of Cartwright's entire functions, i.e.,
entire functions $f$ of exponential type with the property
\begin\eq\label{Cart}
\int_{-\infty}^\infty{\log^+|f(x)|\over 1+x^2}dx<\infty,
\end\eq

by $D$ denote the class of entire functions $f$ of exponential
type with the property: for each $h_m\in\R$ there exists a
subsequence $h'_m$ and a function $\tilde f\not\equiv 0$ such that
the sequence $f(z+h'_m)$ converges uniformly to a function $\tilde
f$ on every compact subset of $\C$.
\medskip

It is easy to see that $D\subset B\subset C$. The classes $A$ and
$B$ are well known and prove useful in various areas of complex
and functional analysis (cf. \cite{L1}, \cite{L}, \cite{K}), class
$D$ is an extension of the well known class of Sine-Type functions
(cf.\cite{L}); it inherits some properties of Sine-Type functions.

Here we will consider sequences $\{a_k\}$ of complex numbers
without finite limit points; note that each $\{a_k\}$ may appear
in the sequence with a finite multiplicity. We say that a sequence
$\{a_k\}$ is {\it a zero set for an entire function} $f$ if
$f(a)\neq 0$ for $a\not\in\{a_k\}$ and
$f(a)=f'(a)=\dots=f^{(p-1)}(a)=0$, $f^{(p)}(a)\neq 0$, where $p$
is a multiplicity of $a$ in $\{a_k\}$.

The zero set of any function $f\in C$ satisfies the conditions
\begin\eq\label{A}
 \sum_{a_k\neq 0}|\Im a_k^{-1}|<\infty,
\end\eq
\begin\eq\label{Lind}
\exists \lim_{R\to\infty}\sum_{0<|a_k|<R} a_k^{-1},
\end\eq
and for all $0<\a\le\pi/2$
\begin\eq\label{reg1}
\lim_{R\to\infty}R^{-1}\card\{a_k:\,0<|a_k|<R,\,|\arg
a_k|\le\a\}=d/2\pi,
\end\eq
\begin\eq\label{reg2}
\lim_{R\to\infty}R^{-1}\card\{a_k:\,0<|a_k|<R,\,|\arg
a_k-\pi|\le\a\}=d/2\pi,
\end\eq
where $d$ is the width of the indicator diagram of $f$
(cf.\cite{L}, p.127). Then every function $f\in C$ has the form
\begin\eq\label{form}
f(z)=Az^se^{i\l z}\lim_{R\to\infty}\prod_{0<|a_k|<R}(1-z/a_k)
\end\eq
with $A\in\C,\,\l\in\R$ (cf.\cite{L}, p.130). However there are
zero sets satisfying (\ref{A})--(\ref{reg2}) such that every
entire function with this zero set does not belong to this class
\footnote{For example, take $a_1=-e^2$ and reals
$a_k<a_1,\;k=2,\dots,$ such that the number
$n(r)=\card\{a_k:\,a_k\ge -r\}$ is the integer part of $r/\log^2r$
for each $r\ge e^2$. Then (\ref{A})--(\ref{reg2}) fulfilled and
the function (\ref{form}) with $s=0$ does not satisfy (\ref{Cart})
because for $x>e^2$ we have
$$
\log|f(x)|=\int_1^\infty\log(1+x/t)dn(t)=\int_1^\infty{xn(t)dt\over
t(x+t)}>1+\int_x^\infty{xdt\over 2t\log^2t}=1+{x\over 2\log x};
$$
the idea belongs to prof. A.\,Grishin.}.

 B.\,Levin in \cite{L1}, Appendix VI,
investigated zero sets $\{a_k\}, k\in\Z$, for functions of class
$B$ with the additional property $\sup_k |a_k-hk|<\infty$ with
some $h\in\R$.  In the present paper we obtain a complete
description of zero sets for the classes $B$, $C$, $D$.

In what follows we put $n(c,t)=\card\{a_k:\,|a_k-c|\le t\}$ for
$c\in\C$. Assume also that the sequence $\{a_k\}$ satisfies the
following conditions
\begin\eq\label{O}
n(0,t)=O(t), \quad  t\to\infty,
\end\eq
and
\begin\eq\label{o}
n(0,t+1)-n(0,t)=o(t), \quad t\to\infty.
\end\eq
It is clear that (\ref{O}) and (\ref{o}) are invariant with
respect to the change $n(0,t)$ to $n(c,t)$.
\begin{Th}\label{1}
A sequence $\{a_k\}\subset\C$ is the zero set of some function
$f\in C$ if and only if it satisfies (\ref{Lind}), (\ref{O}),
(\ref{o}), and
\begin\eq\label{C}
\int_{-\infty}^\infty\left[\int_0^\infty
[n(b,t)-n(x,t)]t^{-1}dt\right]^+{dx\over 1+x^2}<\infty,
\end\eq
for some point $b\in\R\setminus\{a_k\}$.
\end{Th}

\begin{Th}\label{2}
A sequence $\{a_k\}\subset\C$ is the zero set of some function
$f\in B$ if and only if it satisfies (\ref{Lind}), (\ref{O}),
(\ref{o}), and
\begin\eq\label{B}
\sup_{x\in\R}\int_0^\infty [n(b,t)-n(x,t)]t^{-1}dt<\infty
\end\eq
for some point $b\in\R\setminus\{a_k\}$ \footnote {Theorems 1 and
2 were announced (with my permission) in survey \cite{Kh}, p.45.}.

\medskip

\end{Th}

\begin{Th}\label{3}
A sequence $\{a_k\}\subset\C$ is the zero set of some function
$f\in D$ if and only if it satisfies (\ref{Lind}), (\ref{O}),
(\ref{o}), and
\begin\eq\label{D}
\sup_{x\in\R}\left|\int_1^\infty
[n(0,t)-n(x,t)]t^{-1}dt\right|<\infty.
\end\eq
\end{Th}
The proofs of Theorems 1--3 are based on the following lemma:

\begin{Le}. Let a sequence $\{a_k\}\subset\C\setminus\{0\}$
satisfy (\ref{Lind}), (\ref{O}), and (\ref{o}). Then
\begin{\eq}\label{g}
g(z)=\lim_{R\to\infty}\prod_{|a_k|<R}(1-z/a_k)
\end{\eq}
is a well-defined entire function of a finite exponential type,
and for all $z\in\C$
\begin\eq\label{log}
\log|g(z)|=\int_0^\infty [n(0,t)-n(z,t)]t^{-1}dt.
\end\eq
Moreover, if $z_0$ is a point from the sequence $\{a_k\}$ with
multiplicity $l=n(z_0,0)>0$, then
\begin\eq\label{log0}
\log\left|{g^{(l)}(z_0)\over l!}\right|=\int_1^\infty
{n(0,t)-n(z_0,t)\over t}dt+ \int_0^1 {n(0,t)-n(z_0,t)+l\over t}dt.
\end\eq
\end{Le}

{\bf Proof of Lemma}. It follows from (\ref{O}) that the integral
$\int_0^\infty n(0,t)t^{-3}$ is finite, and so
\begin\eq\label{sq}
\sum |a_k|^{-2}=\int_0^\infty t^{-2}dn(0,t)<\infty.
\end\eq
Let $K$ be a fixed disc in the complex plane. Then
$$
\sum_{R<|a_k|<R'}\log\left(1-{z\over
a_k}\right)-\sum_{R<|a_k|<R'}{z\over a_k}\to 0
$$
as $R, R'\to\infty$ uniformly on $z\in K$.  Hence, uniformly on
$z\in K$ we have
\begin\eq\label{un}
\lim_{R\to\infty}\sum_{|a_k|>R}\log\left(1-{z\over a_k}\right)=0,
\end\eq
and the function $g$ is well-defined. Conditions (\ref{Lind}) and
(\ref{O}) now yield that $g$ (\ref{g}) has a finite exponential
type of growth.

Next, for $z=a_k$ equality (\ref{log}) is trivial. Fix
$z\not\in\{a_k\}$. Take $\e>0$ such that $n(0,\e)=n(z,\e)=0$. Then
we have for any $R<\infty$
\begin\eq\label{sum1}
\sum_{|a_k|\le R}(\log|z-a_k|-\log|a_k|)=\sum_{\e\le|a_k-z|\le R}
\log|z-a_k| - \sum_{\e\le|a_k|\le R} \log|a_k| \\
+\sum_{|a_k|\le R, |z-a_k|>R} \log|z-a_k| -\sum_{|z-a_k|\le R,
|a_k|>R} \log|z-a_k|.\nonumber
\end\eq
Furthermore,
$$
\sum_{\e\le|a_k-z|\le R}\log|z-a_k|= n(z,R)\log R- \int_\e^R
n(z,t)t^{-1}dt,
$$
and
$$
(\log R) [n(z,R)-n(0,R)]+ \sum_{|a_k|\le R, |z-a_k|>R} \log|z-a_k|
-\sum_{|z-a_k|\le R, |a_k|>R} \log|z-a_k|
$$
\begin{\eq}\label{sum3}
 =\sum_{|a_k|\le R, |z-a_k|>R}
\log{|z-a_k|\over R} +\sum_{|z-a_k|\le R, |a_k|>R} \log{R\over
|z-a_k|}.
\end{\eq}
Note that if $|z-a_k|>R$ and $|a_k|\le R$, then $|a_k|>R-|z|$ and
$|z-a_k|\le |z|+R$. Also, if $|a_k|>R$ and $|z-a_k|\le R$, then
$|a_k|\le R+|z|$ and $|z-a_k|>R-|z|$. Hence for $R\to\infty$ the
right-hand side of (\ref{sum3}) does not exceed
\begin\eq\label{OO}
O(|z|/R)[n(0,R)-n(0,R-|z|)]+O(|z|/R)[n(0,R+|z|)-n(0,R)].
\end\eq
Taking into account (\ref{o}), we see that (\ref{OO}) tends to
zero for fixed $z$ as $R\to\infty$. Therefore  assertion
(\ref{log}) follows from (\ref{sum1}) and (\ref{un}). In order to
get (\ref{log0}), apply (\ref{log}) to the function
$g(z)(z-z_0)^{-l}z_0^l$. Lemma is proved.

{\bf Proof of Theorems 1--3}. It follows from (\ref{sq}) that the
limits as $R\to\infty$ of the sums $\sum_{|a_k-c|<R}(a_k-c)^{-1}$
and $\sum_{|a_k-c|<R}(a_k)^{-1}$ exist simultaneously, and the
last sum differs from the sum $\sum_{|a_k|<R}(a_k)^{-1}$ no more
than $n(0,R+|c|)-n(0,R-|c|)=o(R)$ terms, where modulus of each
term does not exceed $(R-|c|)^{-1}$. Therefore (\ref{Lind}) is
invariant with respect to the shift of the origin. Further, using
Lemma 1 for the function
$\lim_{R\to\infty}\prod_{|a_k-c|<R}(1-z/(a_k-c))$, we see that
$\left[\int_0^\infty [n(c,t)-n(x,t)]t^{-1}dt\right]^+=O(|x|)$ as
$|x|\to\infty$, therefore (\ref{C}) is invariant with respect to
the shift on a real $b$ as well. Hence with no loss of generality
we can assume $0\not\in\{a_k\}$ and $b=0$ in (\ref{C}), (\ref{B}).

It follows immediately from Lemma 1 that under conditions of
Theorems 1 or 2 the function (\ref{g}) belongs to class $C$ or
$B$, respectively. Let $\{a_k\}$ satisfy the conditions of Theorem
3. By (\ref{log}), (\ref{log0}), and Jensen formula
$$
{1\over 2\pi}\int_0^{2\pi}\log|g(z+e^{i\theta})|d\theta=
\log\left|{g^{(l)}(z)\over l!}\right|+\int_0^1 {n(z,t)-n(z,0)\over
t}dt
$$
with $l=n(z,0)$, we get
\begin\eq\label{log2}
{1\over
2\pi}\int_0^{2\pi}\log|g(z+e^{i\theta})|d\theta=\int_1^\infty
{n(0,t)-n(z,t)\over t}dt+\int_0^1 {n(0,t)\over t}dt.
\end\eq
Hence it follows from (\ref{D}) and the inequality
$$
\log|g(z)|\le {1\over
2\pi}\int_0^{2\pi}\log|g(z+e^{i\theta})|d\theta,
$$
that $g$ is bounded on $\R$. The Pfragmen--Lindelof Theorem
implies that $g(z)$ is bounded on every horizontal strip of a
finite width. Therefore for each $t_m\in\R$ there exists a
subsequence $t'_m$ such that $f(z+t'_m)$ converge uniformly to a
function $h(z)$ on every compact subset of $\C$. If $h(z)\equiv 0$
on $\C$, then $\int_0^{2\pi}\log|g(t_m'+e^{i\theta})|d\theta\to
-\infty$. By (\ref{log2}), this contradicts (\ref{D}), and so
$g\in D$.

To prove necessity, note that by (\ref{form}) every function $f\in
C$ with $f(0)\neq 0$ has the form
$$
f(z)=Ae^{i\l z}g(z)
$$
with  $g$ defined in (\ref{g}). The conditions (\ref{O}) and
(\ref{o}) follow immediately from (\ref{reg1}) and (\ref{reg2}).
Using Lemma 1, we see that the zero set of $f$ satisfies (\ref{C})
and, when $f\in B$, (\ref{B}). For the case $f\in D$ the functions
$f$ and $g$ are bounded on $\R$ and, by the Pfragmen--Lindelof
Theorem, on every horizontal strip of a finite width. Hence, if
(\ref{D}) is false, then from (\ref{log2})
$$
{1\over 2\pi}\int_0^{2\pi}\log|g(t_l+e^{i\theta})|d\theta\to
-\infty
$$
for some sequence of reals $t_j$. Taking into account properties
of averages of subharmonic functions, we have
$$
\int_{|u+iv-t_j|<1}\log|g(u+iv)|du\,dv\to -\infty.
$$
Since the function $\log|g|$ is bounded from above on the strip
$|\Im z|<2$, we get
$$
\int_{|u+iv-t_j-s|<2}\log|g(u+iv)|du\,dv\to -\infty
$$
uniformly in $s\in [-1,1]$. Using properties of averages of
subharmonic functions again, we obtain
\begin\eq\label{inf}
\sup_{-1\le s\le 1}\log|g(t_j+s)|\to -\infty.
\end\eq
On the other hand, $f(z+t'_j)\to h(z)\not\equiv 0$ for some
subsequence $t'_j$ uniformly on compact subsets of $\C$, hence
(\ref{inf}) is impossible. Theorems 1--3 are completely proved.
\medskip

{\bf Remark 1}. It can be shown that all these theorems are valid
with the function $\tilde n(c,t)=\card\{a_k:\,|\Re(a_k-c)|\le t,
|\Im(a_k-c)|\le t\}$ instead of $n(c,t)$ as well.

{\bf Remark 2}. The additional condition
$$
\limsup_{y\to\pm\infty}\int_0^\infty
[n(b,t)-n(iy,t)](|y|t)^{-1}dt\le\sigma
$$
in Theorems 1--3 gives the complete description for zero sets of
corresponding subclasses of functions with an exponential type at
most $\sigma$.
\medskip

 {\bf Example}. Let $\a(t)$ be a strictly
increasing concave function on $\R^+$ with the properties
$\a(0)=0$, $1\le\a'(t)\le 1+O(t^{-1})$ as $t\to\infty$, $\a(t)\ge
t+1+\e$ for large $t$ with some $\e>0$. Consider a sequence of
reals $\{a_k\},\, k\in\Z\setminus\{0\}$, such that $\a(a_k)=k$ and
$\{a_{-k}\}=\{-a_k\},\, k\in\N$. Note that
$$
n(0,t)-n(x,t)=2\E[\a(t)]+\E[\a(x-t)]-\E[\a(x+t)]
$$
for $x>t$, and
$$
n(0,t)-n(x,t)=2\E[\a(t)]-E[\a(x-t)]-\E[\a(x+t)]
$$
for $0<x<t$; here $\E[x]$ is the integer part of $x$. Since
$a_{n+1}-a_n=(\a'(\tilde t))^{-1}$ with some $\tilde t\in
(a_n,\,a_{n+1})$, we see that for $n\ge n_0$
$$
0\le 1-(a_{n+1}-a_n)\le(a_{n+1}-a_n)^{-1}-1\le
\a'(a_n)-1<C/a_n<1/2,
$$
and
$$
\left|\int_{a_n}^{a_{n+1}} (\a(t)-\E[\a(t)]-1/2)dt\right|\le
\left|\int_{a_n}^{a_{n+1}} \left(\int_{a_n}^t
\a'(u)du-1/2\right)dt\right|\le C'/a_n.
$$
Therefore we have for arbitrary $N<\infty$
$$
\left|\int_{a_{n_0}}^{a_N}(\a(t)-\E[\a(t)]-1/2)dt\right|\le C'\log
a_N.
$$
Now it follows easily that all the integrals
$$
\int_1^x {\a(t)-\E[\a(t)]-1/2\over t}dt,\quad \int_1^x {\a(x\pm
t)-\E[\a(x\pm t)]-1/2\over t}dt,
$$
$$
\int_x^\infty {\a(t\pm x)-\E[\a(t\pm x)]-1/2\over t}dt,
$$
are bounded uniformly in $x\to\infty$. Hence, up to some uniformly
bounded term, integral $\int_0^\infty [n(0,t)-n(x,t)]t^{-1}dt$
equals
\begin\eq\label{int}
\int_1^x {2\E[\a(t)]-2t\over
t}dt+\int_1^{x-1}{\a(x-t)-\a(x+t)+2t\over t}dt\nonumber\\
+\int_x^\infty {2\a(t)-\a(t-x)-\a(x+t)\over t}dt.
\end\eq
The function $\a(t)$ is concave, therefore the last integral in
(\ref{int}) is positive. Further, we have for $1<t<x-1$
$$
\a(x-t)-\a(x+t)+2t=\int_{-t}^t 1-\a'(x+u)du\ge -C\log{x+t\over
x-t}.
$$
Now it is easy to check that the second integral in (\ref{int}) is
uniformly bounded from below. Evidently, the first integral in
(\ref{int}) is unbounded from above as $x\to\infty$. Thus the
sequence $\{a_k\}$ is the zero set of no function $f\in B$.

\bigskip

Department of Mathematics, Kharkov National University,

Svobody sq.,4, Kharkov 61077, Ukraine

\end{document}